\documentclass[11pt]{article}
\usepackage[margin=1in]{geometry}
\usepackage{amsmath,amssymb,amsthm,graphicx,hyperref,booktabs}
\hypersetup{hidelinks}
\newtheorem{proposition}{Proposition}[section]
\newtheorem{theorem}[proposition]{Theorem}
\newtheorem{corollary}[proposition]{Corollary}
\newtheorem{lemma}[proposition]{Lemma}
\theoremstyle{remark}
\newtheorem{remark}[proposition]{Remark}
\newcommand{\norm}[1]{\left\lVert#1\right\rVert}
\title{Forward and Inverse Born Series for Scattering Problems with
Second Harmonic Generation}
\author{Ross McNeill$^{1}$ \qquad Shari Moskow$^{1}$ \qquad
John C. Schotland$^{2}$\\[0.75ex]
\small $^{1}$Department of Mathematics, Drexel University,
Philadelphia, Pennsylvania 19104, USA\\
\small $^{2}$Department of Mathematics and Department of Physics,\\
\small Yale University, New Haven, Connecticut 06511, USA}
\begin{document}
\maketitle
\begin{abstract}
We study the forward and inverse Born series for a scalar model of
second harmonic generation (SHG), in which a quadratic susceptibility couples
the Helmholtz equations for the fields at the fundamental and second harmonic
frequencies and the linear susceptibility is known. The inverse problem is
to recover the quadratic susceptibility from boundary measurements of the two
fields. We write the forward Born
series for this system in terms of multilinear operators and obtain explicit
bounds on these operators, which give an explicit radius of convergence.
Convergence of the inverse Born series for sufficiently small data then
follows from a theorem of Hoskins and Schotland. The model has an exact
scaling property: the data depend on the incident field strength and the
susceptibility only through their product. As a consequence, for exact data
the inverse Born series converges for a susceptibility of any size, provided
the incident field is weak enough. Numerical reconstructions in
two dimensions from multifrequency data show that susceptibilities of order
one are recovered with a weak source, and that the higher order terms of the
series substantially improve on the linear reconstruction. We also present a
numerical study of the full model, in which the linear susceptibility is
unknown as well.
\end{abstract}

\section{Introduction}\label{sec:intro}
There has been considerable recent interest in inverse problems for nonlinear
partial differential equations. For elliptic and hyperbolic equations with
power type nonlinearities, it is known that the nonlinearity can be used
to recover coefficients that would not be determined by the corresponding
linear problem \cite{Kurylev2018,Lassas2021}. Much less is known about
practical reconstruction methods for such problems. In previous work
\cite{DeFilippis2023,DeFilippis2024}, we developed the forward and inverse Born
series for scalar waves with a Kerr (cubic) nonlinearity. The inverse Born
series (IBS) expresses the solution of the inverse problem as an explicitly
computable functional of the measured data. It was introduced for diffuse
waves \cite{Markel2003,Moskow2008}, and its convergence and stability have
been studied for diffuse and scalar waves
\cite{Moskow2008,Kilgore2012,Mahankali2023}, for the Calder\'on problem
\cite{Arridge2012,Abhishek2020}, for problems with internal data
\cite{Bardsley2014}, on graphs \cite{Chung2017}, and in a general setting
using geometric function theory \cite{Hoskins2022}. The extension to nonlinear problems requires the forward
Born series to be written in terms of multilinear operators, and the
combinatorics are different from the linear case.

In this paper we consider second harmonic generation, the simplest
nonlinear optical process, in which a quadratic susceptibility converts light
at frequency $\omega$ into light at frequency $2\omega$ \cite{Boyd2008}.
Second harmonic generation is the basis of SHG microscopy
\cite{Campagnola2003}. Inverse problems for SHG have been studied for
Maxwell's equations \cite{Assylbekov2021}, where uniqueness of the linear and
nonlinear parameters was shown, and for a coupled Helmholtz system with
internal data \cite{Ren2024}. Here we study the scalar SHG model with the 
linear susceptibility known, and without the complex conjugate in the equation for the
fundamental field. This is what we refer to as the reduced model. Omitting the
conjugate makes the fields complex analytic functions of the quadratic
susceptibility, which is what the convergence theorem of \cite{Hoskins2022}
requires. We reconstruct the quadratic susceptibility from boundary
measurements of both harmonics.

Our results are as follows. We write the Born series for the two fields in
terms of multilinear operators and obtain bounds for these operators
(Theorem~\ref{thm:forward}). The proof uses a generating function that turns out to be the same one that
appears for the Kerr nonlinearity \cite{DeFilippis2024}, and the bounds give
an explicit radius of convergence for the forward series. Together with the
convergence theorem of
\cite{Hoskins2022}, they give convergence of the inverse Born series for
sufficiently small data (Corollary~\ref{cor:shg-convergence}).

The reduced model has an exact scaling property
(Proposition~\ref{prop:scaling}): the data for a source $\gamma g$ and
susceptibility $\eta$ are $\gamma$ times the data for the source $g$ and
susceptibility $\gamma\eta$. We therefore can show that, for exact data, the IBS
converges for any fixed susceptibility if the incident field is weak enough
(Corollary~\ref{cor:scaling}). This is the SHG analog of the result of
\cite{DeFilippis2024}. We also observe that the first order data only
involve the second harmonic, and that as a consequence all even order terms of
the IBS vanish and the fundamental frequency data are not used
(Proposition~\ref{rem:parity}). Finally, we present numerical reconstructions in
two dimensions, in which the data are computed on a finer grid than the one
used for the inversion, and in which the higher order terms of the series
substantially improve the reconstruction. We also present a numerical study of
the full SHG model, in which both the linear and the quadratic
susceptibilities are unknown, and the complex conjugate appears.

The paper is organized as follows. In Section~\ref{sec:model} we describe the
model and its scaling property. The forward Born series is studied in
Section~\ref{sec:forward}, and the inverse Born series in
Section~\ref{sec:inverse}. Numerical results are presented in
Section~\ref{sec:numerics}. In Section~\ref{sec:full} we present a numerical
study of the full SHG model.

\section{The reduced SHG model}\label{sec:model}
Let $\Omega\subset\mathbb R^d$, $d=2,3$, be a bounded domain with smooth
boundary, and let $\nu$ denote the outward unit normal on $\partial\Omega$. Let
$k>0$ be the wavenumber at the fundamental frequency, let $g$ be a prescribed
boundary source, and let $\eta\in L^\infty(\Omega)$ be the quadratic
susceptibility. The fields $u_1$ and $u_2$ at the fundamental and second
harmonic frequencies satisfy
\begin{equation}\label{eq:pde}
\begin{aligned}
(\Delta+k^2)u_1&=-8\pi k^2\eta u_1u_2 \quad\text{in }\Omega,&
\partial_\nu u_1-iku_1&=g \quad\text{on }\partial\Omega,\\
(\Delta+4k^2)u_2&=-16\pi k^2\eta u_1^2 \quad\text{in }\Omega,&
\partial_\nu u_2-2iku_2&=0 \quad\text{on }\partial\Omega.
\end{aligned}
\end{equation}
This is the scalar SHG model (see Section~\ref{sec:full}) with known linear
susceptibility, with two simplifications. First, we take the known linear
susceptibility to be zero. This is only for convenience: all of the results
hold for any known linear susceptibility (Remark~\ref{rem:known-linear}).
Second, in the first equation we have written $u_1u_2$ in place of
$\overline{u_1}u_2$. Without the complex conjugate the fields are complex
analytic functions of $\eta$, so that the Born series are power series in
$\eta$ with complex multilinear coefficients; the model with the conjugate
is discussed in Remark~\ref{rem:conjugate}. We use the time dependence
$e^{-i\omega t}$.
The impedance boundary conditions are the simplest absorbing boundary
conditions: they are the outgoing (Sommerfeld) radiation condition
$\partial_ru-iku\to0$ imposed on $\partial\Omega$, and they model waves leaving
the domain. Unlike Neumann or Dirichlet conditions, they give a well-posed
linear problem for every $k>0$. Indeed, if $(\Delta+k^2)w=0$ in $\Omega$ and
$\partial_\nu w-ikw=0$ on $\partial\Omega$, then multiplying by $\overline w$
and integrating by parts gives
$\int_\Omega(|\nabla w|^2-k^2|w|^2)=ik\int_{\partial\Omega}|w|^2$, so $w=0$ on
$\partial\Omega$; then also $\partial_\nu w=0$ there, and $w\equiv0$ by unique
continuation. Thus there are no real values of $k$ at which the problem fails
to be uniquely solvable, in contrast to the Neumann eigenvalues, and we will
see in Section~\ref{sec:numerics} that this is important in computations at
higher frequencies. All of the analysis below holds equally well for Neumann
conditions, provided $k^2$ and $4k^2$ are not Neumann eigenvalues.

For $j=1,2$ we define the resolvent $R_j$ by
\begin{equation}
R_j f=w \quad \text{if and only if} \quad
\begin{cases}
(\Delta+j^2k^2)w=f&\text{in }\Omega,\\
\partial_\nu w-ijkw=0&\text{on }\partial\Omega.
\end{cases} \label{eq:resolvent-pde}
\end{equation}
Equivalently, $(R_jf)(x)=\int_\Omega G_j(x,y)f(y)\,dy$, where $G_j$ is the
corresponding Green's function. Since $G_j(x,\cdot)$ has an integrable
singularity, $R_j$ maps $L^\infty(\Omega)$ boundedly into
$C(\overline\Omega)$, with
\begin{equation}
\norm{R_j}=\norm{R_j}_{L^\infty\to C}=\sup_{x\in\Omega}\int_\Omega|G_j(x,y)|\,dy.
\label{eq:resolvent-bounded}
\end{equation}
The background field $u_{10}$ solves $(\Delta+k^2)u_{10}=0$ in $\Omega$ with
$\partial_\nu u_{10}-iku_{10}=g$ on $\partial\Omega$, and we set
$u_0=(u_{10},0)$. There is no background field at the second harmonic
frequency; $u_2$ is generated entirely by the nonlinearity. It follows from
\eqref{eq:pde} that $u_1=u_{10}-8\pi k^2R_1(\eta u_1u_2)$ and
$u_2=-16\pi k^2R_2(\eta u_1^2)$. We write these integral equations as the
system
\begin{equation}\label{eq:int}
\begin{pmatrix}u_1\\ u_2\end{pmatrix}
=\begin{pmatrix}u_{10}\\ 0\end{pmatrix}
-\begin{pmatrix}0&8\pi k^2R_1(\eta u_1\,\cdot\,)\\[0.5ex]
16\pi k^2R_2(\eta u_1\,\cdot\,)&0\end{pmatrix}
\begin{pmatrix}u_1\\ u_2\end{pmatrix},
\end{equation}
where $R_j(\eta u_1\,\cdot\,)$ denotes the operator $v\mapsto R_j(\eta u_1v)$.
That is,
\begin{equation}\label{eq:int-T}
u=u_0+T_\eta(u),\qquad u=\begin{pmatrix}u_1\\ u_2\end{pmatrix},\qquad
u_0=\begin{pmatrix}u_{10}\\ 0\end{pmatrix},
\end{equation}
where, for $v=(v_1,v_2)$,
\begin{equation}\label{eq:T-def}
T_\eta(v)=-\begin{pmatrix}0&8\pi k^2R_1(\eta v_1\,\cdot\,)\\[0.5ex]
16\pi k^2R_2(\eta v_1\,\cdot\,)&0\end{pmatrix}
\begin{pmatrix}v_1\\ v_2\end{pmatrix}
=\begin{pmatrix}-8\pi k^2R_1(\eta v_1v_2)\\[0.5ex]
-16\pi k^2R_2(\eta v_1^2)\end{pmatrix}.
\end{equation}
The matrix has no diagonal entries: the fundamental field is driven only by
the second harmonic field, and the second harmonic field only by the
fundamental field. The map $T_\eta$ is quadratic in $v$ and linear in $\eta$.
We work in the space $X=C(\overline\Omega)^2$ with the norm
\begin{equation}\label{eq:Xnorm}
\norm{(v_1,v_2)}_X=\norm{v_1}_\infty+\norm{v_2}_\infty,
\end{equation}
and we define
\begin{equation}
\mu_1=8\pi k^2\norm{R_1},\qquad
\mu_2=16\pi k^2\norm{R_2},\qquad \mu=\max(\mu_1,\mu_2).
\label{eq:resolventconstants}
\end{equation}

\begin{proposition}\label{prop:wellposed}
Let $u_{10}$ be the background field (the solution of \eqref{eq:pde} with
$\eta=0$), let $u_0=(u_{10},0)$, so that
$\norm{u_0}=\norm{u_0}_X=\norm{u_{10}}_\infty$, and let $\mu$ be the constant in
\eqref{eq:resolventconstants}. If $\eta\in L^\infty(\Omega)$ satisfies
\begin{equation}\label{eq:wellposed-condition}
4\mu\norm{\eta}_\infty\norm{u_0}<1,
\end{equation}
then the integral equations \eqref{eq:int} for $(u_1,u_2)$ have a unique
solution in the ball
\[
\overline B=\{v\in X:\norm{v-u_0}_X\leq\norm{u_0}\}
\]
of $X=C(\overline\Omega)^2$.
\end{proposition}
\begin{proof}
By \eqref{eq:int-T}, a solution of \eqref{eq:int} is a fixed point of the map
$v\mapsto u_0+T_\eta(v)$. If
$v\in\overline B$ then $\norm{v_1}+\norm{v_2}\leq 2\norm{u_0}$, and so
$\norm{T_\eta(v)}_X\leq
\mu\norm{\eta}_\infty\norm{v_1}(\norm{v_1}+\norm{v_2})
\leq 4\mu\norm{\eta}_\infty \norm{u_0}^2<\norm{u_0}$. For $v,w\in\overline B$, using
$v_1v_2-w_1w_2=(v_1-w_1)v_2+w_1(v_2-w_2)$ and
$v_1^2-w_1^2=(v_1+w_1)(v_1-w_1)$, we obtain
\[
\begin{aligned}
\norm{T_\eta(v)-T_\eta(w)}_X
&\leq\mu\norm{\eta}_\infty\Big[\big(\norm{v_1}+\norm{v_2}+\norm{w_1}\big)
\norm{v_1-w_1}+\norm{w_1}\norm{v_2-w_2}\Big]\\
&\leq 4\mu\norm{\eta}_\infty \norm{u_0}\,\norm{v-w}_X .
\end{aligned}
\]
Thus $v\mapsto u_0+T_\eta(v)$ is a contraction of $\overline B$ into itself.
\end{proof}

The data consist of measurements for a finite number of acquisitions
$\ell=1,\ldots,L$. We label acquisition $\ell$ to have source $g_\ell$ and a wavenumber
$k_\ell$ (several sources may share the same wavenumber). The measurements are the
boundary values of the scattered fields at both frequencies,
$\phi_\ell=\big((u_1-u_{10})|_{\partial\Omega},\,u_2|_{\partial\Omega}\big)
\in C(\partial\Omega)^2$, where $u_1$, $u_2$ and $u_{10}$ are the fields for
that source and wavenumber. The data are
\begin{equation}
\phi=(\phi_1,\ldots,\phi_L)\in Y=\big(C(\partial\Omega)^2\big)^L,\qquad
\norm{\phi}_Y=\max_{1\leq\ell\leq L}\big(\norm{\phi_{\ell,1}}_\infty
+\norm{\phi_{\ell,2}}_\infty\big),
\label{eq:data-definition}
\end{equation}
where $\phi_\ell=(\phi_{\ell,1},\phi_{\ell,2})$. The fields $u_{10}$ depend on the
acquisition, and the resolvents $R_j$ and the constants $\mu_j$ depend on its
wavenumber. To keep the notation simple, we write the analysis for a single
acquisition, with source $g$ and wavenumber $k$. For the full data set, we
use the same symbols for the largest of these constants,
\begin{equation}\label{eq:maxconstants}
\norm{u_0}=\max_\ell\norm{u_{10}^{(\ell)}}_\infty,\qquad
\mu_1=\max_\ell8\pi k_\ell^2\norm{R_1^{(\ell)}},\qquad
\mu_2=\max_\ell16\pi k_\ell^2\norm{R_2^{(\ell)}},
\end{equation}
where the superscript indicates the acquisition. This is consistent with the
maximum norm on $Y$, and all of the estimates for a single acquisition remain
true with these larger constants. When
we need to indicate the dependence on the sources and the susceptibility, we
write $\phi=\phi[g,\eta]$, with $g=(g_1,\ldots,g_L)$. The following scaling property is the key to
several results below.

\begin{proposition}\label{prop:scaling}
Let $\gamma>0$ and $\eta\in L^\infty(\Omega)$. If $(u_1,u_2)$ solves the
reduced SHG system \eqref{eq:pde} with source $g$ and susceptibility
$\gamma\eta$, then $(\gamma u_1,\gamma u_2)$ solves \eqref{eq:pde} with source
$\gamma g$ and susceptibility $\eta$. In particular, for the data
\eqref{eq:data-definition} with sources $g=(g_1,\ldots,g_L)$, if all of the
sources are multiplied by $\gamma$, then
\begin{equation}
\phi[\gamma g,\eta]=\gamma\,\phi[g,\gamma\eta]. \label{eq:scaling-identity}
\end{equation}
\end{proposition}
\begin{proof}
The boundary conditions are linear, and
\[
(\Delta+k^2)(\gamma u_1)=-8\pi k^2(\gamma\eta)\gamma u_1u_2
=-8\pi k^2\eta(\gamma u_1)(\gamma u_2).
\]
The second equation can be checked similarly.
\end{proof}

Thus the strength of the nonlinear effect is determined by the product of
the incident field and the susceptibility, as one expects from the physics of
SHG \cite{Boyd2008}. We will refer to $\gamma\norm{\eta}_\infty$ as the
effective contrast.

\begin{remark}[Known linear susceptibility]\label{rem:known-linear}
We have taken the linear susceptibility to be zero, but what matters is that
it is known. Suppose that a real linear susceptibility
$\eta_1\in L^\infty(\Omega)$ is given, and replace the operators $\Delta+k^2$
and $\Delta+4k^2$ in \eqref{eq:pde} by $\Delta+k^2(1+4\pi\eta_1)$ and
$\Delta+4k^2(1+4\pi\eta_1)$. This is the model \eqref{eq:fullpde} of
Section~\ref{sec:full} with $\eta_1$ known and $\eta=\eta_2$ unknown (without
the complex conjugate). The uniqueness argument above applies without change,
since the additional term $4\pi k^2\int_\Omega\eta_1|w|^2$ is real, so the
linear problems are again well posed for every $k>0$. We then define the
resolvents $R_j$ in \eqref{eq:resolvent-pde} with the new operators, and we
let $u_{10}$ be the field in the known medium, that is, the solution with
$\eta=0$. The integral equations \eqref{eq:int} have exactly the same form,
and the rest of the analysis only uses this form and the norms $\norm{R_j}$.
Hence Propositions~\ref{prop:wellposed} and~\ref{prop:scaling}, and the
results of Sections~\ref{sec:forward} and~\ref{sec:shg-ibs}, hold with the
same proofs, with the constants $\mu_j$ and $\norm{u_0}$ computed for the known medium.
In particular, the parity of the fields, the fact that the IBS does not use
the fundamental frequency data, and the convergence for weak incident fields
do not depend on the linear susceptibility being zero, only on its being
known.
\end{remark}

\section{The forward Born series}\label{sec:forward}
We now write the Born series in terms of multilinear operators, following
\cite{DeFilippis2023,DeFilippis2024}. If
$L:L^\infty(\Omega)^n\to Y$ (or $C(\overline\Omega)$, or $X$) is $n$-linear, we use the norm
\[
\norm{L}=\sup\big\{\norm{L(\eta_1,\ldots,\eta_n)}:\norm{\eta_j}_\infty\leq1,\
j=1,\ldots,n\big\}.
\]
We use the notation of \cite{DeFilippis2023}. For $v\in C(\overline\Omega)$
and $\zeta\in L^\infty(\Omega)$, let
\begin{equation}\label{eq:bj}
b_j(v,\zeta)=-8\pi jk^2\,R_j(\zeta v),\qquad j=1,2,
\end{equation}
where $R_j$ is the resolvent defined in \eqref{eq:resolvent-pde}.
These operators have tensor counterparts. Given a multilinear operator
$T_l:L^\infty(\Omega)^l\to C(\overline\Omega)$ of order $l$, we define the
multilinear operator $B_jT_l$ of order $l+1$ by
\begin{equation}\label{eq:liftings}
B_jT_l(\zeta_1,\ldots,\zeta_l,\zeta_{l+1})
=b_j\big(T_l(\zeta_1,\ldots,\zeta_l),\zeta_{l+1}\big),\qquad j=1,2.
\end{equation}
Given multilinear operators $T_l$ and $T_m$ of orders $l$ and $m$, their tensor
product is the multilinear operator of order $l+m$ given by
\begin{equation}\label{eq:tensor}
T_l\otimes T_m(\zeta_1,\ldots,\zeta_l,\zeta_{l+1},\ldots,\zeta_{l+m})
=T_l(\zeta_1,\ldots,\zeta_l)\,T_m(\zeta_{l+1},\ldots,\zeta_{l+m}).
\end{equation}
An operator of order zero is a fixed function. In this notation, the integral
equations \eqref{eq:int} read
\[
u_1=u_{10}+B_1\,u_1\otimes u_2,\qquad u_2=B_2\,u_1\otimes u_1,
\]
where the fields are viewed as sums of multilinear operators in $\eta$. By
\eqref{eq:resolvent-bounded} and \eqref{eq:resolventconstants},
\begin{equation}\label{eq:Bbound}
\norm{T_l\otimes T_m}\leq\norm{T_l}\norm{T_m},\qquad
\norm{B_jT_l}\leq\mu_j\norm{T_l}.
\end{equation}
The operators $U_{j,n}:L^\infty(\Omega)^n\to C(\overline\Omega)$ are defined by
$U_{1,0}=u_{10}$, $U_{2,0}=0$, and, for $n\geq0$,
\begin{equation}\label{eq:forwardops}
U_{1,n+1}
=\sum_{p+p'=n}B_1\,U_{1,p}\otimes U_{2,p'},\qquad
U_{2,n+1}
=\sum_{p+p'=n}B_2\,U_{1,p}\otimes U_{1,p'},
\end{equation}
where the sums are over $p,p'\geq0$. The forward operators $K_n:L^\infty(\Omega)^n\to Y$ are
the $n$-linear operators
\begin{equation}
K_n(\eta_1,\ldots,\eta_n)
=\big(U_{1,n}(\eta_1,\ldots,\eta_n)|_{\partial\Omega},\,
       U_{2,n}(\eta_1,\ldots,\eta_n)|_{\partial\Omega}\big), \label{eq:Kn}
\end{equation}
and the Born series for the data is
\begin{equation}
\phi=\sum_{n=1}^\infty K_n(\eta,\ldots,\eta). \label{eq:dataseries}
\end{equation}
Here, as explained in Section~\ref{sec:model}, the operators are written for a
single acquisition. For the full data set,
$K_n=(K_n^{(1)},\ldots,K_n^{(L)})$, where $K_n^{(\ell)}$ is defined as above with
the source $g_\ell$ and wavenumber $k_\ell$, and since $Y$ has the maximum norm,
$\norm{K_n}=\max_\ell\norm{K_n^{(\ell)}}$.

\begin{lemma}[Parity]\label{lem:parity}
Let $U_{1,n}$ and $U_{2,n}$ be the operators for the fundamental and second
harmonic fields defined by \eqref{eq:forwardops}, with $U_{1,0}=u_{10}$ and
$U_{2,0}=0$. Then $U_{1,n}=0$ for $n$ odd, and $U_{2,n}=0$ for $n$ even. In
particular, for each acquisition the first order operator \eqref{eq:Kn} is
$K_1(\eta)=\big(0,\,-16\pi k^2R_2(\eta u_{10}^2)|_{\partial\Omega}\big)$, where
$R_2$ is the resolvent \eqref{eq:resolvent-pde} at wavenumber $2k$, so it only
has second harmonic data.
\end{lemma}
\begin{proof}
We use induction on $n$. For $n=0$, $U_{2,0}=0$. Suppose the statement holds
for all orders up to $n$. A term $B_1\,U_{1,p}\otimes U_{2,p'}$ with
$p+p'=n$ can only be nonzero if $p$ is even and $p'$ is odd, so that $n$ is
odd. Hence $U_{1,n+1}=0$ if $n+1$ is odd. Similarly, a term
$B_2\,U_{1,p}\otimes U_{1,p'}$ can only be nonzero if $p$ and $p'$ are both
even, so that $n$ is even, and hence $U_{2,n+1}=0$ if $n+1$ is even. The
formula for $K_1$ follows from \eqref{eq:forwardops} with $n=0$.
\end{proof}

That is, the fundamental field only contains even orders in $\eta$, and the
second harmonic field only contains odd orders. By
Proposition~\ref{prop:scaling}, if the source $g$ is replaced by $\gamma g$,
then $K_n$ is replaced by $\gamma^{n+1}K_n$.

\subsection{Bounds on the forward operators}
We now bound the forward operators, which gives an explicit radius of
convergence for the Born series.

\begin{theorem}\label{thm:forward}
Let $U_{1,n}$, $U_{2,n}$ and $K_n$ be the forward operators of the reduced SHG
problem, defined by \eqref{eq:forwardops} and \eqref{eq:Kn}, with the norms of
multilinear operators defined at the beginning of this section and the data
space $Y$ of \eqref{eq:data-definition}. Let $\norm{u_0}$ be the largest
value of $\norm{u_{10}}_\infty$ over the acquisitions, where $u_{10}$ is the
background field, and let $\mu_1$ and $\mu_2$ be the largest values of
$8\pi k^2\norm{R_1}$ and $16\pi k^2\norm{R_2}$, as in
\eqref{eq:maxconstants}. Let
\begin{equation}\label{eq:forwardconstants}
\mu_B=\frac{3\sqrt3}{2}\,\norm{u_0}\sqrt{\mu_1\mu_2},\qquad
\nu_B=\frac{\sqrt3}{2}\,\norm{u_0}\Big(\sqrt3+\sqrt{\frac{\mu_2}{\mu_1}}\Big),
\qquad r_F=\frac1{\mu_B}.
\end{equation}
Then
\begin{equation}
\norm{K_n}\leq\nu_B\mu_B^n,\qquad n\geq1. \label{eq:geometricKn}
\end{equation}
In particular, if
\[
\norm{\eta}_\infty\leq r_F,
\]
the Born series $u_j=\sum_nU_{j,n}(\eta,\ldots,\eta)$ for the fields converges absolutely in
$C(\overline\Omega)$ for every acquisition, its sum $(u_1,u_2)$ is a solution
of the integral equations \eqref{eq:int}, and the Born series
\eqref{eq:dataseries} for the data converges absolutely in $Y$.
\end{theorem}
\begin{proof}
\emph{Step 1: a recursion for bounds.} Define numbers $p_n,q_n\geq0$ by
$p_0=\norm{u_0}$, $q_0=0$, and
\begin{equation}\label{eq:pqrec}
p_{n+1}=\mu_1\sum_{i+j=n}p_iq_j,\qquad
q_{n+1}=\mu_2\sum_{i+j=n}p_ip_j,
\qquad n\geq0.
\end{equation}
These numbers do not depend on the acquisition. Fix an acquisition. We show by
induction that $\norm{U_{1,n}}\leq p_n$ and $\norm{U_{2,n}}\leq q_n$. This
holds for $n=0$, since $\norm{U_{1,0}}=\norm{u_{10}}_\infty\leq\norm{u_0}$ and
$U_{2,0}=0$.
If it holds up to order $n$, then by \eqref{eq:forwardops} and
\eqref{eq:Bbound},
\[
\norm{U_{1,n+1}}\leq\mu_1\sum_{p+p'=n}\norm{U_{1,p}}\norm{U_{2,p'}}
\leq\mu_1\sum_{p+p'=n}p_pq_{p'}=p_{n+1},
\]
and in the same way $\norm{U_{2,n+1}}\leq q_{n+1}$. Here we used that
$\mu_1$ and $\mu_2$ are maxima over the acquisitions
\eqref{eq:maxconstants}.

\emph{Step 2: generating functions.} As in \cite{DeFilippis2023,DeFilippis2024},
we bound the operators through generating functions. We work with formal power
series in a variable $t$. Let $P(t)=\sum_{n\geq0}p_nt^n$ and $Q(t)=\sum_{n\geq0}q_nt^n$.
The sums in \eqref{eq:pqrec} are the coefficients of the power series for
the products $PQ$ and $P^2$:
\[
P(t)Q(t)=\sum_{n\geq0}\Big(\sum_{i+j=n}p_iq_j\Big)t^n,\qquad
P(t)^2=\sum_{n\geq0}\Big(\sum_{i+j=n}p_ip_j\Big)t^n .
\]
Multiplying the first relation in \eqref{eq:pqrec} by $t^{n+1}$ and summing
over $n\geq0$ therefore gives
\[
P-\norm{u_0}=\sum_{n\geq0}p_{n+1}t^{n+1}
=\mu_1t\sum_{n\geq0}\Big(\sum_{i+j=n}p_iq_j\Big)t^n=\mu_1tPQ,
\]
and in the same way, since $q_0=0$, the second relation gives $Q=\mu_2tP^2$.
Substituting the second identity into the first, we obtain
\begin{equation}
P=\norm{u_0}+\mu_1tPQ,\qquad Q=\mu_2tP^2,
\qquad\text{and hence}\qquad P=\norm{u_0}+\mu_1\mu_2t^2P^3. \label{eq:PQcubic}
\end{equation}

\emph{Step 3: the series $W$.} Let $W(x)=\sum_{m\geq0}\tau_mx^m$ be the formal
power series solution of
\begin{equation}\label{eq:ternary}
W=1+xW^3.
\end{equation}
It exists and is unique, since comparing coefficients of $x^m$ determines
$\tau_m$ from $\tau_0,\ldots,\tau_{m-1}$, starting from $\tau_0=1$, and all of
the coefficients are positive. Dividing the last identity in
\eqref{eq:PQcubic} by $\norm{u_0}$ shows that $P/\norm{u_0}$ satisfies
\eqref{eq:ternary} with $x=\mu_1\mu_2\norm{u_0}^2t^2$, so by uniqueness
\begin{equation}\label{eq:PQW}
P(t)=\norm{u_0}\,W(\mu_1\mu_2\norm{u_0}^2t^2),\qquad
Q(t)=\mu_2tP^2=\mu_2\norm{u_0}^2t\,W(\mu_1\mu_2\norm{u_0}^2t^2)^2 .
\end{equation}

\emph{Step 4: the radius of convergence.} The same cubic equation arises for
the Kerr nonlinearity: the generating function $P$ in the proof of
Proposition~2 of \cite{DeFilippis2024} is $\nu_0W(\nu_0^2x)$ in our notation.
It is shown there that the series for $W$ converges for $0\leq x<4/27$, with
$W(x)<\frac32$. Since the coefficients are positive, the partial sums at
$x=4/27$ are then also bounded by $\frac32$, so the series converges at
$x=4/27$, and its sum is a positive root of $\frac4{27}W^3-W+1
=\frac4{27}(W-\frac32)^2(W+3)$, that is, $W(4/27)=\frac32$. Now let
$\rho=1/\mu_B=2/(3\sqrt3\,\norm{u_0}\sqrt{\mu_1\mu_2})$. Then
$\mu_1\mu_2\norm{u_0}^2\rho^2=4/27$, so by \eqref{eq:PQW} the series $P(t)$ and $Q(t)$
converge for $|t|\leq\rho$, and
\[
P(\rho)=\frac{3\norm{u_0}}2,\qquad
Q(\rho)=\mu_2\norm{u_0}^2\rho\cdot\frac94=\frac{\sqrt3\,\norm{u_0}}{2}\sqrt{\frac{\mu_2}{\mu_1}},
\qquad\text{so that}\qquad P(\rho)+Q(\rho)=\nu_B.
\]

\emph{Step 5: conclusion.} Since all coefficients are nonnegative,
$(p_n+q_n)\rho^n\leq P(\rho)+Q(\rho)=\nu_B$, that is,
$p_n+q_n\leq\nu_B\mu_B^n$. The data are the boundary values of the fields, and the
maximum of a continuous function over $\partial\Omega$ is at most its maximum
over $\overline\Omega$. So for every acquisition
\[
\norm{K_n^{(\ell)}}\leq\norm{U_{1,n}}+\norm{U_{2,n}}\leq p_n+q_n
\leq\nu_B\mu_B^n .
\]
Taking the maximum over
$\ell$ gives \eqref{eq:geometricKn}, and the Born series converges for every
acquisition when $\norm{\eta}_\infty\leq r_F$. Inserting the series into the
right-hand side of \eqref{eq:int} and collecting the terms of each order,
which is permitted since the series converge absolutely, gives the recursion
\eqref{eq:forwardops}. Hence the sum satisfies \eqref{eq:int}.
\end{proof}

The constants $\mu_B$ and $\nu_B$ are proportional to $\norm{u_0}$. In particular the radius
$r_F=1/\mu_B$ is inversely proportional to $\norm{u_0}$, so for any given
$\eta\in L^\infty(\Omega)$ the forward Born series converges if the incident
field is weak enough. This is the scaling of
Proposition~\ref{prop:scaling} seen in the bounds. The corresponding
statement for the inverse series also involves $\norm{K_1^+}$, and is given
in Corollary~\ref{cor:scaling}.

Thus, for $\norm{\eta}_\infty\leq r_F$, the sum of the Born series is a
solution of the integral equations \eqref{eq:int}. If $\eta$ also satisfies
the condition \eqref{eq:wellposed-condition}, one checks from the bounds
in the proof of Theorem~\ref{thm:forward} that the sum lies in the ball $\overline B$ of
Proposition~\ref{prop:wellposed}, so it is the solution given by that
proposition. Indeed, the iterates of the fixed-point iteration
$v\mapsto u_0+T_\eta(v)$ of Proposition~\ref{prop:wellposed}, started from the
background field, agree with the Born series through the order of the
iterate. Because of the nonlinearity they also contain some terms of higher
order, so they are not the partial sums of the series.

\section{The inverse Born series}\label{sec:inverse}\label{sec:shg-ibs}
We now turn to the inverse problem. Let $Z=L^\infty(\Omega)$, let $Y$ be the
data space $\big(C(\partial\Omega)^2\big)^L$ of \eqref{eq:data-definition},
and let $K_n$ be the forward operators of Section~\ref{sec:forward}, so that
the forward map is
\begin{equation}
F(\eta)=K_1\eta+\sum_{n=2}^{\infty}K_n(\eta,\ldots,\eta).
\label{eq:forward-splitting}
\end{equation}
The data are $\phi=F(\eta)$, possibly with noise. Let $K_1^+:Y\to Z$ be a
bounded linear operator, which is a regularized inverse of $K_1$, for example
a truncated SVD or a Tikhonov regularized inverse.

\subsection{Definition of the series}
The inverse Born series is $\widetilde\eta=\sum_{m\geq1}\mathcal K_m(\phi)$,
where
\begin{equation}\label{eq:ibs}
\mathcal K_1(\phi)=K_1^+\phi,\qquad
\mathcal K_m(\phi)=-K_1^+
\sum_{n=2}^m\ \sum_{i_1+\cdots+i_n=m}
K_n\big(\mathcal K_{i_1}(\phi),\ldots,\mathcal K_{i_n}(\phi)\big),\quad m\geq2,
\end{equation}
and the inner sum is over $i_1,\ldots,i_n\geq1$. The inverse Born series was
first analyzed in \cite{Moskow2008}, where it is written as a recursion for
the operators $\mathcal K_m$. Here we use the form of the series introduced
in \cite{Hoskins2022} (see also \cite{DeFilippis2023,DeFilippis2024}), in
which each term is obtained from the previous terms
$\mathcal K_1(\phi),\ldots,\mathcal K_{m-1}(\phi)$ themselves, so that these
can be stored and reused.

When $K_1^+$ is a truncated SVD, we take the singular system
$K_1v_j=s_jw_j$ of $K_1$ as a real-linear operator, since the susceptibility
is real, with the $L^2(\Omega)$ inner product on the susceptibilities and the
real inner product $\operatorname{Re}\langle\cdot,\cdot\rangle$ of
$L^2(\partial\Omega)^{2L}$ on the data. If $J$ is the set of retained indices,
then
\begin{equation}\label{eq:real-tsvd}
K_1^+y=\sum_{j\in J}s_j^{-1}\operatorname{Re}\langle y,w_j\rangle\,v_j .
\end{equation}
The retained singular vectors $v_j$ are real, bounded functions, so $K_1^+$
is a bounded operator of finite rank from $Y$ to $L^\infty(\Omega)$. Its range
is $E=\operatorname{span}\{v_j:j\in J\}$, and $K_1^+K_1$ is the
$L^2$-orthogonal projection onto $E$. Every term of the inverse Born series
lies in $E$.

\subsection{The fundamental frequency data}
The parity of the fields (Lemma~\ref{lem:parity}) has a consequence for the
inverse series, which is special to this problem. Write
$\phi=(\phi_1,\phi_2)$, where $\phi_1$ collects the data at the fundamental
frequency and $\phi_2$ the second harmonic data, for all acquisitions. We say
that $K_1^+$ \emph{does not use the fundamental data} if
\begin{equation}\label{eq:no-fund}
K_1^+(\phi_1,0)=0\qquad\text{for all }\phi_1 .
\end{equation}
This is a natural property, since by Lemma~\ref{lem:parity} the first order
data $K_1\eta$ have no fundamental component. It holds for the truncated
SVD, for Tikhonov regularization, and for any other regularized inverse of the
form $K_1^+=f(K_1^*K_1)K_1^*$, when the adjoint is computed with an inner
product in which the two components of the data are orthogonal (as is the case
after discretization). Indeed, the range of $K_1$ is then orthogonal to all
data of the form $(\phi_1,0)$, and $K_1^*$ vanishes on the orthogonal
complement of the range of $K_1$. We will assume \eqref{eq:no-fund} from now
on.

\begin{proposition}[Even order terms]\label{rem:parity}
Consider the reduced SHG problem, with forward operators $K_n$ defined by
\eqref{eq:forwardops}--\eqref{eq:Kn} and inverse Born series terms
$\mathcal K_m(\phi)$ defined by \eqref{eq:ibs}, and let
$K_1^+:Y\to L^\infty(\Omega)$ be a bounded linear operator that does not use
the fundamental data, that is, \eqref{eq:no-fund} holds. Then:
\begin{enumerate}
\item[(i)] $K_1^+K_{2m}=0$ for all $m\geq1$;
\item[(ii)] $\mathcal K_m(\phi)=0$ for all even $m$, and the inverse Born
series only depends on the second harmonic data $\phi_2$.
\end{enumerate}
\end{proposition}
\begin{proof}
(i) By Lemma~\ref{lem:parity}, $K_{2m}$ only has a fundamental component, so
$K_1^+K_{2m}=0$ by \eqref{eq:no-fund}.

(ii) We use induction on $m$. For $m=2$, $\mathcal K_2(\phi)=-K_1^+K_2(\ldots)=0$
by (i). Suppose $\mathcal K_j(\phi)=0$ for all even $j<m$, and let $m$ be even.
In \eqref{eq:ibs}, a term $K_n(\mathcal K_{i_1}(\phi),\ldots,\mathcal K_{i_n}(\phi))$
can only be nonzero if all $i_j$ are odd, and then $n\equiv m\equiv0$ (mod 2).
By (i), $K_1^+K_n=0$ for $n$ even, so $\mathcal K_m(\phi)=0$. Finally,
$\mathcal K_1(\phi)=K_1^+(0,\phi_2)$ by \eqref{eq:no-fund}, and every term of
\eqref{eq:ibs} depends on $\phi$ only through $\mathcal K_1(\phi)$ and the
earlier terms, so the series only depends on $\phi_2$.
\end{proof}

In particular, the measured fundamental data are not used, and the nonlinear
terms of the inverse series start at cubic order.

\subsection{Convergence}
Convergence of the inverse Born series follows from the bounds on the forward
operators in Theorem~\ref{thm:forward} and the following theorem, which we state in our notation, as in
\cite{DeFilippis2023,DeFilippis2024}: a bound $\nu\mu^n$ here is a bound
$\nu'\mu^{n-1}$ in \cite{Hoskins2022} with $\nu'=\nu\mu$. An important step in
\cite{Hoskins2022} was the rewriting of the series in the form \eqref{eq:ibs},
which removed the smallness assumption on $\norm{K_1^+}$ required in
\cite{Moskow2008}. The proof in \cite{Hoskins2022} uses holomorphic maps on
complex Banach spaces, so we state the theorem for complex spaces and
complex-multilinear operators.

\begin{theorem}[\cite{Hoskins2022}]\label{thm:HS}
Let $Z$ and $Y$ be complex Banach spaces, let the forward map be
$F(\eta)=\sum_{n\geq1}K_n(\eta,\ldots,\eta)$ as in
\eqref{eq:forward-splitting}, with $K_n:Z^n\to Y$ bounded and complex
$n$-linear, let $K_1^+:Y\to Z$ be bounded and complex-linear, and let
$\phi\in Y$. Suppose that
$\norm{K_n}\leq\nu\mu^n$ for all $n\geq1$, and let
\[
C=\max\{2,\norm{K_1^+}\nu\mu\},\qquad
\rho=\frac{\sqrt{16C^2+1}-4C}{2\mu}.
\]
If
\[
\norm{K_1^+\phi}<\rho,
\]
then the inverse Born series
$\sum_m\mathcal K_m(\phi)$, with terms defined by \eqref{eq:ibs}, converges.
\end{theorem}

\begin{corollary}\label{cor:shg-convergence}
Consider the reduced SHG problem, with forward operators $K_n$ defined by
\eqref{eq:forwardops}--\eqref{eq:Kn}, data $\phi$ in the space $Y$ of
\eqref{eq:data-definition}, and the constants $\nu_B$ and $\mu_B$ of
Theorem~\ref{thm:forward}. Let $K_1^+:Y\to L^\infty(\Omega)$ be a bounded
linear operator, or the truncated SVD \eqref{eq:real-tsvd}, and let
\[
C=\max\{2,\norm{K_1^+}\nu_B\mu_B\},\qquad
\rho=\frac{\sqrt{16C^2+1}-4C}{2\mu_B}.
\]
If
\[
\norm{K_1^+\phi}_\infty<\rho,
\]
then the inverse Born series \eqref{eq:ibs} converges in $L^\infty(\Omega)$.
\end{corollary}
\begin{proof}
By Theorem~\ref{thm:forward}, $\norm{K_n}\leq\nu_B\mu_B^n$ for all $n\geq1$,
so Theorem~\ref{thm:HS} applies with $\nu=\nu_B$ and $\mu=\mu_B$. For the
truncated SVD \eqref{eq:real-tsvd}, which is only real-linear, see
Remark~\ref{rem:real-tsvd}.
\end{proof}

\begin{remark}[Real susceptibilities]\label{rem:real-tsvd}
Theorem~\ref{thm:HS} requires a complex-linear operator $K_1^+$, while
the truncated SVD \eqref{eq:real-tsvd} is only real-linear on $Y$, because of
the real part. The theorem can still be applied, as follows. Let
$\widehat Y=Y\times Y$ with the norm
$\max(\norm{y_1}_Y,\norm{y_2}_Y)$, and define
\begin{align*}
\widehat K_n(\eta_1,\ldots,\eta_n)
&=\Big(K_n(\eta_1,\ldots,\eta_n),\
\overline{K_n(\overline{\eta_1},\ldots,\overline{\eta_n})}\Big),\\
\widehat K_1^+(y_1,y_2)&=\frac12\sum_{j\in J}s_j^{-1}
\Big(\langle y_1,w_j\rangle+\langle y_2,\overline{w_j}\rangle\Big)v_j .
\end{align*}
The operators $\widehat K_n$ are complex $n$-linear, with
$\norm{\widehat K_n}=\norm{K_n}$, and $\widehat K_1^+$ is complex-linear. For
real $\eta_1,\ldots,\eta_n$ we have
$\widehat K_n(\eta_1,\ldots,\eta_n)=(y,\overline y)$ with
$y=K_n(\eta_1,\ldots,\eta_n)$, and
$\widehat K_1^+(y,\overline y)=K_1^+y$, since
$\langle\overline y,\overline{w_j}\rangle=\overline{\langle y,w_j\rangle}$. It
follows by induction that the inverse Born series for the operators
$\widehat K_n$, the regularized inverse $\widehat K_1^+$ and the data
$(\phi,\overline\phi)$ has the same terms as the series \eqref{eq:ibs}, all of
which are real. Finally, $\norm{\widehat K_1^+}\leq\norm{K_1^+}$. Indeed, let
$K_{1,\mathbb C}^+y=\sum_{j\in J}s_j^{-1}\langle y,w_j\rangle v_j$, so that
$K_1^+y=\operatorname{Re}K_{1,\mathbb C}^+y$. For $x\in\Omega$ and $y\in Y$,
choose $\theta$ with
$|(K_{1,\mathbb C}^+y)(x)|=e^{-i\theta}(K_{1,\mathbb C}^+y)(x)
=\big(K_1^+(e^{-i\theta}y)\big)(x)\leq\norm{K_1^+}\norm{y}_Y$. Hence
$\norm{K_{1,\mathbb C}^+}\leq\norm{K_1^+}$, and the same bound holds for the
second term in $\widehat K_1^+$, which is
$\overline{K_{1,\mathbb C}^+\overline{y_2}}$. Thus Theorem~\ref{thm:HS},
applied to the operators $\widehat K_n$ and $\widehat K_1^+$, gives the
statement of Corollary~\ref{cor:shg-convergence} for the truncated SVD.
\end{remark}

Theorem~\ref{thm:HS} gives convergence, but it does not identify the limit of
the series. Estimates of the difference between the limit and the true
coefficient, in terms of $\norm{(I-K_1^+K_1)\eta}$, are given in
\cite{Moskow2008,Hoskins2022}. For the truncated SVD, all partial sums lie in
the space $E$, so the best that can be expected is the $L^2$-orthogonal
projection of the true susceptibility onto $E$, and this is what we compare
with in Section~\ref{sec:numerics}.

\begin{remark}[The complex conjugate]\label{rem:conjugate}
In the scalar SHG model the first equation in \eqref{eq:pde} has
$\overline{u_1}u_2$ in place of $u_1u_2$ (see Section~\ref{sec:full}). With
the conjugate and $\eta$ real, the fields are no longer complex analytic in
$\eta$, but they are real-analytic, and the forward operators are
real-multilinear. Some of the results above carry over. The parity of
Lemma~\ref{lem:parity} is easy to see directly in this case: if $(u_1,u_2)$
is the solution for $\eta$, then $(u_1,-u_2)$ is the solution for $-\eta$, so
$u_1$ is even and $u_2$ is odd in $\eta$. The scaling of
Proposition~\ref{prop:scaling} holds as well, and the estimates of
Section~\ref{sec:forward} only involve absolute values of the fields, so we
expect the bounds on the forward operators to be unchanged, although we do
not write out the details. However, the proof of Theorem~\ref{thm:HS} uses
holomorphic maps on complex Banach spaces, so the theorem applies to
complex-multilinear forward operators, and
Corollary~\ref{cor:shg-convergence} does not apply as stated to the model
with the conjugate. We expect that this is a technical restriction, and that
the convergence result holds for that model as well.
\end{remark}

\subsection{Weak incident fields}
In \cite{DeFilippis2024} we showed for Kerr and more general polynomial
nonlinearities, with known linear term, that the region of convergence of the
IBS grows as the incident field is made weaker. For the reduced SHG problem
this follows directly from Proposition~\ref{prop:scaling}. Let
$\mathcal K_m^{\gamma}$ denote the IBS operators for the source $\gamma g$, and
let $K_1^+$ be a truncated SVD with a cutoff relative to the largest singular
value.

\emph{Step 1: scaling of the operators.} By Proposition~\ref{prop:scaling},
the forward operators for the source $\gamma g$ are $\gamma^{n+1}K_n$. In
particular the first order operator is $\gamma^2K_1$, whose singular values
are $\gamma^2$ times those of $K_1$, with the same singular vectors. With a
relative cutoff the same singular vectors are retained, so the regularized
inverse is $\gamma^{-2}K_1^+$ and the space $E$ does not depend on $\gamma$.

\emph{Step 2: scaling of the IBS.} We show by induction that
\begin{equation}
\mathcal K_m^{\gamma}(\phi)=\gamma^{-1}\mathcal K_m^{1}(\phi/\gamma).
\label{eq:ibs-scaling}
\end{equation}
For $m=1$, $\mathcal K_1^\gamma(\phi)=\gamma^{-2}K_1^+\phi
=\gamma^{-1}K_1^+(\phi/\gamma)$. If \eqref{eq:ibs-scaling} holds for all orders
less than $m$, then by \eqref{eq:ibs} and multilinearity,
\[
\mathcal K_m^\gamma(\phi)
=-\gamma^{-2}K_1^+\sum_{n,\,i_1+\cdots+i_n=m}\gamma^{n+1}\gamma^{-n}
K_n\big(\mathcal K^1_{i_1}(\phi/\gamma),\ldots,\mathcal K^1_{i_n}(\phi/\gamma)\big)
=\gamma^{-1}\mathcal K^1_m(\phi/\gamma).
\]
Combined with \eqref{eq:scaling-identity}, this says that reconstructing $\eta$
with the source $\gamma g$ is exactly the same as reconstructing $\gamma\eta$
with the source $g$, followed by division by $\gamma$.

\begin{corollary}\label{cor:scaling}
Consider the reduced SHG problem with sources $g=(g_1,\ldots,g_L)$ and data
$\phi[g,\eta]$ as in \eqref{eq:data-definition}. Let $K_1^+$ be the truncated
SVD of the first order operator $K_1$ for the sources $g$, with a cutoff
relative to the largest singular value. There is $\rho_0>0$, depending only
on $g$ and on the cutoff, such that if $\eta_*\in L^\infty(\Omega)$ and
$\gamma>0$ satisfy
\[
\gamma\norm{\eta_*}_\infty<\rho_0 ,
\]
then the inverse Born series \eqref{eq:ibs} for the sources $\gamma g$ and
the exact data $\phi[\gamma g,\eta_*]$ converges. In particular, for every
$\eta_*\in L^\infty(\Omega)$ the series converges if the sources are weak
enough.
\end{corollary}
\begin{proof}
By \eqref{eq:scaling-identity} and \eqref{eq:ibs-scaling},
\[
\sum_m\mathcal K_m^\gamma\big(\phi[\gamma g,\eta_*]\big)
=\gamma^{-1}\sum_m\mathcal K_m^1\big(\phi[g,\gamma\eta_*]\big),
\]
so it is enough to consider the IBS for the sources $g$ and the data
$\phi[g,\zeta]$, where $\zeta=\gamma\eta_*$, and to show that it converges when
$s=\norm{\zeta}_\infty$ is small. By Theorem~\ref{thm:forward}, for
$s<r_F=1/\mu_B$ we have $\phi[g,\zeta]=\sum_{n\geq1}K_n(\zeta,\ldots,\zeta)$
and
\[
\norm{K_1^+\phi[g,\zeta]}_\infty
\leq\norm{K_1^+}\sum_{n\geq1}\nu_B\mu_B^ns^n
=\norm{K_1^+}\frac{\nu_B\mu_Bs}{1-\mu_Bs},
\]
which tends to zero with $s$. Hence
$\norm{K_1^+\phi[g,\zeta]}_\infty<\rho$ for $s$ small enough, where
$\rho$ is the radius of Corollary~\ref{cor:shg-convergence} for the
sources $g$, and the series converges by that corollary.
\end{proof}

Thus, for exact data, a large susceptibility is not in itself an obstruction
to convergence of the IBS, since the incident field can be made weaker. There
is, however, an issue with noise. If the incident field is multiplied by
$\gamma$, the second harmonic data, which are what the series uses, are
multiplied by $\gamma^2$ to leading order. If the noise has a fixed size
$\delta$ that does not depend on the incident field, then its size relative
to the data is of order $\delta/\gamma^2$. So $\gamma$ cannot be taken
arbitrarily small: it has to be small enough for the series to converge, and
large enough that the second harmonic signal, of order $\gamma^2$, remains
well above the noise level $\delta$.

\section{Numerical results}\label{sec:numerics}
\subsection{Setup}
We take $\Omega=(0,1)^2$ and discretize \eqref{eq:pde} on a uniform grid with
the five-point Laplacian. The boundary conditions are imposed by eliminating
ghost points, and $\eta$ is represented by its values at the interior nodes.
The nonlinear forward problem is solved by iterating \eqref{eq:int}, with
Newton's method as a fallback, to a relative tolerance of $10^{-12}$. The
terms of the IBS are those defined by \eqref{eq:ibs}; by
Proposition~\ref{rem:parity}, only the terms of odd order are nonzero.
Computing them only requires solves with the background operators at $k$ and
$2k$, for which we use sparse LU factorizations.

To avoid an inverse crime, the data are computed on a $128\times128$ grid,
while the inversion uses a $64\times64$ grid, so that there are
$62^2=3844$ unknowns. The two grids have $127$ and $63$ intervals in each
direction, which are relatively prime, so no interior nodes coincide. The
susceptibility is given by a formula and sampled separately on the two grids,
and the boundary values on the fine grid are interpolated to the boundary
nodes of the coarse grid. We use $12$ equally spaced wavenumbers from
$k=0.75$ to $k=9.5$. For each wavenumber there are five sources on each side
of the square, centered at $0.1,0.3,0.5,0.7,0.9$; each source is a Gaussian
with standard deviation $0.1$ and integral one, multiplied by the source
strength $\gamma$. This gives $240$ source--frequency pairs. By
Proposition~\ref{rem:parity}, only the second harmonic data are used, and we
add complex Gaussian noise of relative size $10^{-3}$ to these data.

The matrix $K_1$ has $240\times252$ rows and is not formed. The row of $K_1$
for a boundary node $x_b$ and a source with background field $u_{10}$ is the
product of $u_{10}^2$ with the solution of the adjoint problem with a point
source at $x_b$. These adjoint solutions only depend on the frequency, so
products with $K_1^*$ are computed with one adjoint solve per boundary node and
frequency. Since $\eta$ is real, we regard $K_1$ as a real-linear map, as
described in Section~\ref{sec:shg-ibs}. Its Gram matrix is then
$\operatorname{Re}(K_1^*K_1)$, and the truncated SVD is obtained from the
eigenvalue decomposition of this real symmetric matrix, so that the
reconstructions and the space $E$ are real. We use a cutoff of $0.03$ relative
to the largest singular value, which keeps $277$ of the $3844$ singular
values. As a measure of the
discretization error, the first order data computed on the two grids differ
by about one percent for the first example below, and by a few percent for
the second, where the susceptibility is discontinuous.

We note that with Neumann boundary conditions we were not able to use
wavenumbers above about $k=6$. The closed domain has no loss, and near the
eigenvalues of the Laplacian the solution is very sensitive to small changes,
including the change from one grid to the other, so that for larger
wavenumbers the first order data on the two grids differed substantially,
even with fine grids. The impedance conditions remove this problem, and the
difference between the grids then grows only slowly with $k$.

In the tables, ``IBS $N$'' is the partial sum of the series through order $N$,
and ``projection'' is the orthogonal projection $P_E\eta$ of the true
susceptibility onto the retained singular vectors. Since every partial sum
lies in $E$, no partial sum can have a smaller $L^2$ error than the
projection. We compute the
relative $L^2$ error $\norm{\eta_N-\eta}/\norm{\eta}$, and also the error
relative to the projection, $\norm{\eta_N-P_E\eta}/\norm{P_E\eta}$, which
mainly measures the nonlinear error, separately from the truncation error.

\subsection{A Gaussian inclusion}
In the first example, the susceptibility is
\begin{equation}
\eta(x,y)=\exp\left(
-\frac{(x-0.371)^2+(y-0.433)^2}{2(0.07)^2}\right).
\label{eq:gaussian-target}
\end{equation}
which has maximum value one. By Proposition~\ref{prop:scaling}, the data for
the source $\gamma g$ and susceptibility $\eta$ are $\gamma$ times the data
for the source $g$ and susceptibility $\gamma\eta$, so the strength of the
nonlinear effects is determined by the effective contrast
$\gamma\norm{\eta}_\infty$, which here is equal to $\gamma$. We therefore
fix $\eta$ and vary the source strength, as in
Corollary~\ref{cor:scaling}, and we recover a susceptibility of order one by using a weak source. 

\begin{table}[ht]
\centering
\begin{tabular}{lrrrrrrr}
\toprule
 & Proj. & IBS 1 & IBS 3 & IBS 5 & IBS 7 & IBS 9 & IBS 11\\
\midrule
Error, $\gamma=0.35$ & $2.0\%$ & $4.6\%$ & $2.4\%$ & $2.3\%$ & $2.3\%$ & $2.3\%$ & $2.3\%$\\
Error, $\gamma=0.5$ & $2.0\%$ & $9.2\%$ & $4.6\%$ & $4.0\%$ & $3.6\%$ & $3.2\%$ & $3.2\%$\\
\quad rel.\ to projection & -- & $9.0\%$ & $4.2\%$ & $3.5\%$ & $3.0\%$ & $2.5\%$ & $2.5\%$\\
Error, $\gamma=0.6$ & $2.0\%$ & $12.0\%$ & $11.7\%$ & $12.2\%$ & $12.2\%$ & $12.0\%$ & $14.0\%$\\
\bottomrule
\end{tabular}
\caption{Gaussian inclusion: relative $L^2$ error of the projection and of the
partial sums of the IBS for three source strengths. For $\gamma=0.5$ we also
give the error relative to the projection.}
\label{tab:gaussian}
\end{table}

\begin{figure}[ht]
\centering
\includegraphics[width=\textwidth]{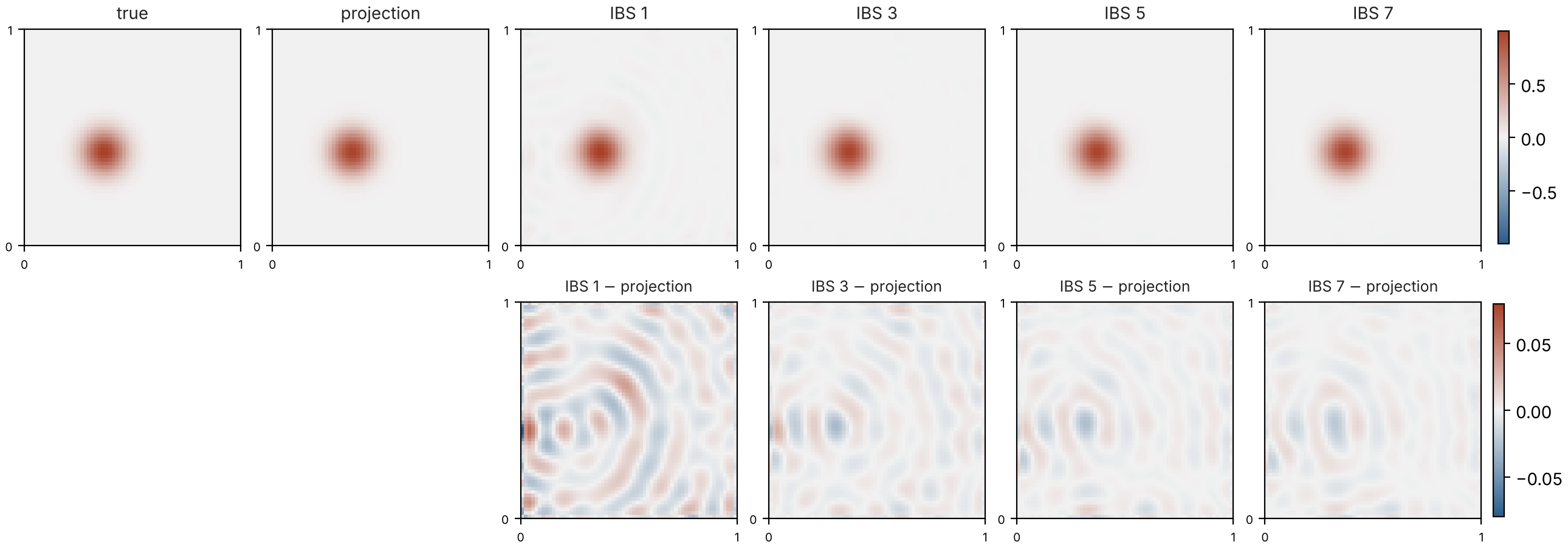}
\caption{Reconstruction of the Gaussian \eqref{eq:gaussian-target} with
$\gamma=0.5$, for which the nonlinear part is about one fifth of the data.
Top: the true susceptibility, its projection onto the retained
singular vectors, and partial sums of the IBS, on the same color scale.
Bottom: the difference between each partial sum and the projection, on a
finer color scale. The data are computed on a $128\times128$ grid and the
inversion uses a $64\times64$ grid.}
\label{fig:gaussian}
\end{figure}

\begin{figure}[ht]
\centering
\includegraphics[width=0.55\textwidth]{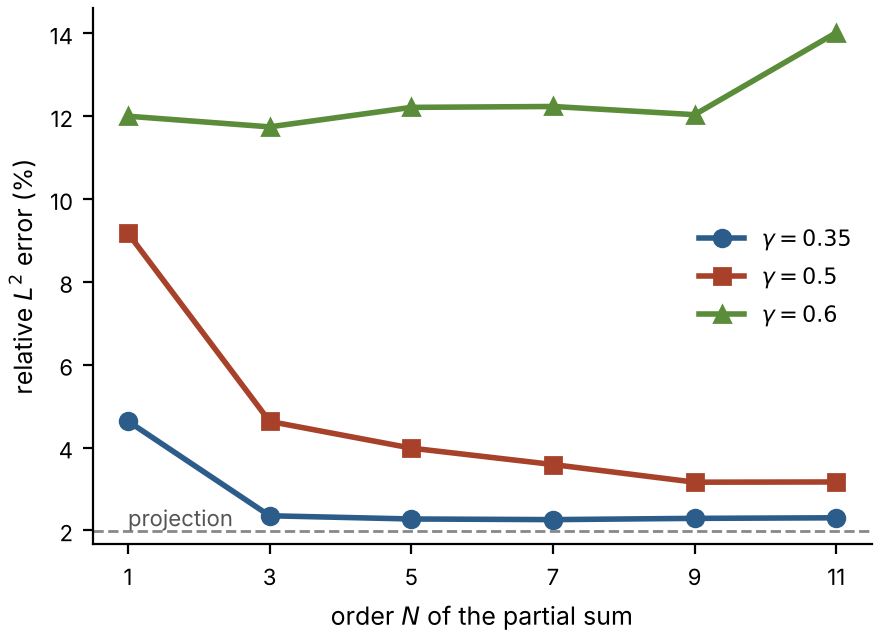}
\caption{Relative $L^2$ error of the partial sums of the IBS for the Gaussian
inclusion, for three source strengths. The dashed line is the error of the
projection.}
\label{fig:gaussian-orders}
\end{figure}

The results are shown in Table~\ref{tab:gaussian} and
Figures~\ref{fig:gaussian} and~\ref{fig:gaussian-orders}. For $\gamma=0.35$
the error of the first term is more than twice that of the projection, and
the series converges quickly to a reconstruction that is nearly as accurate
as the projection; the small remaining difference is consistent with the
discretization error and the noise. For $\gamma=0.5$ the higher order terms
make a large improvement: the error is reduced to about a third of that of
the first term by the partial sum of order nine, and it is then stable. The
error of the first term consists of oscillating artifacts spread over the
domain, which are visible in the bottom row of Figure~\ref{fig:gaussian}, and
the higher order terms remove most of them. The peak value and width of the
reconstruction are close to those of the projection for all partial sums. The
convergence is slower than for $\gamma=0.35$. For $\gamma=0.6$, the higher
order terms no longer improve the reconstruction: the error stays near that
of the first term and then increases, so the series diverges. This is the
behavior expected from Corollary~\ref{cor:shg-convergence} and
Corollary~\ref{cor:scaling}: the size of $\norm{K_1^+\phi}$ relative to the
radius of convergence is determined by the effective contrast, and as it
increases the terms decay more slowly, until the series diverges.

\subsection{A piecewise constant susceptibility}
In the second example the susceptibility is $1$ on a disk of radius $0.15$
centered at $(0.34,0.42)$, $0.6$ on the rectangle
$[0.59,0.81]\times[0.52,0.80]$, and zero elsewhere. The source strength is
$\gamma=0.3$, for which the nonlinear part is again about one fifth of the
data.

\begin{table}[ht]
\centering
\begin{tabular}{lrrrrrrr}
\toprule
 & Proj. & IBS 1 & IBS 3 & IBS 5 & IBS 7 & IBS 9 & IBS 11\\
\midrule
Error & $32.2\%$ & $34.0\%$ & $33.3\%$ & $33.3\%$ & $33.3\%$ & $33.4\%$ & $33.9\%$\\
Rel.\ to projection & -- & $11.6\%$ & $9.1\%$ & $9.2\%$ & $8.9\%$ & $9.5\%$ & $11.3\%$\\
\bottomrule
\end{tabular}
\caption{Piecewise constant susceptibility, $\gamma=0.3$: relative $L^2$ error
of the projection and of the partial sums of the IBS.}
\label{tab:piecewise}
\end{table}

\begin{figure}[ht]
\centering
\includegraphics[width=\textwidth]{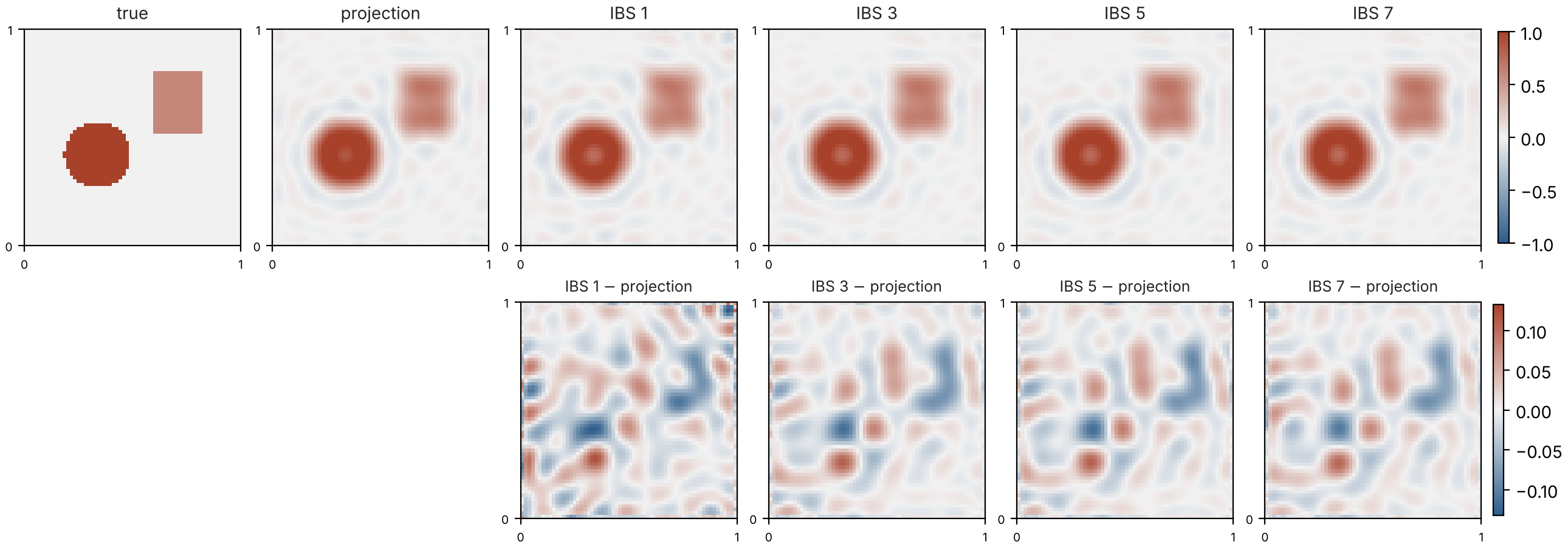}
\caption{Reconstruction of a piecewise constant susceptibility with
$\gamma=0.3$. Top: the true susceptibility, its projection, and partial sums of
the IBS. Bottom: the difference between each partial sum and the projection,
on a finer color scale.}
\label{fig:piecewise}
\end{figure}

Here the error is dominated by the truncation: the projection smooths the
edges of the disk and the rectangle, and its error is large (see
Table~\ref{tab:piecewise} and Figure~\ref{fig:piecewise}). Relative to the
projection, the nonlinear error decreases modestly from the first term to the
partial sums of orders three to seven, after which it increases. The
improvement is smaller than for the Gaussian. One reason is that the part of
$\eta$ outside of $E$, which is large here, still contributes to the
nonlinear terms in the data, and the series cannot account for it, since it
only reconstructs the part in $E$. The model error is also larger for this
discontinuous susceptibility. Thus the nonlinear correction is limited by the
resolution of the linear inversion.

\section{The full SHG model: a numerical study}\label{sec:full}
The scalar SHG model with both linear and quadratic susceptibilities is
\begin{equation}\label{eq:fullpde}
\begin{aligned}
(\Delta+k^2)u_1&=-4\pi k^2\big(\eta_1u_1+2\eta_2u_2\overline{u_1}\big),\\
(\Delta+4k^2)u_2&=-16\pi k^2\big(\eta_1u_2+\eta_2u_1^2\big),
\end{aligned}
\end{equation}
with the same boundary conditions as in \eqref{eq:pde}, where $\eta_1$ is the
linear susceptibility and $\eta_2$ the quadratic one. Both are now unknown.
We do not develop the theory for this model here, but the method of the
previous sections can be carried over, and in this section we describe how,
and present a numerical study.

\subsection{The inverse Born series for the full model}
\emph{Step 1: the Born series.} We take $\eta_1$ and $\eta_2$ to be real. Because
of the complex conjugate in \eqref{eq:fullpde}, the fields are not complex
analytic functions of $(\eta_1,\eta_2)$, but they are real-analytic. As in
Section~\ref{sec:forward}, we expand the fields in terms that are
homogeneous of degree $n$ in $(\eta_1,\eta_2)$, starting from $u_1=u_{10}$
and $u_2=0$ at order zero. The boundary values of the terms of degree $n$
define forward operators $K_n$ that are real-multilinear in
$(\eta_1,\eta_2)$, and the inverse Born series is defined by the same
formula \eqref{eq:ibs}.

\emph{Step 2: the first order operator.} To first order in
$(\eta_1,\eta_2)$, the fields are $u_1\approx u_{10}-4\pi k^2R_1(\eta_1u_{10})$
and $u_2\approx-16\pi k^2R_2(\eta_2u_{10}^2)$. Thus $K_1$ is block diagonal: to first
order, $\eta_1$ is determined by the fundamental data (through the usual
linear Born approximation at wavenumber $k$), and $\eta_2$ by the second
harmonic data, as in the reduced model. We use a separate truncated SVD for
each block, computed from the real Gram matrix
$\operatorname{Re}(K_1^*K_1)$ of that block, as in Section~\ref{sec:numerics}.

\emph{Step 3: differences from the reduced model.} Since $\eta_1$ is now
unknown (compare Remark~\ref{rem:known-linear}), the fundamental data are
used, and there is no parity, so all orders contribute. Moreover, the
terms involving only $\eta_1$ have the same degree in the incident field, so
the scaling of Proposition~\ref{prop:scaling} does not apply to them: making
the source weaker reduces the effective quadratic contrast
$\gamma\norm{\eta_2}_\infty$, but not the linear contrast $\norm{\eta_1}_\infty$.
We therefore expect convergence to require a bound on $\eta_1$, as for scalar
waves \cite{Kilgore2012} and for the Kerr problem with linear and nonlinear
terms \cite{DeFilippis2023}, even if a weak incident field still allows a
large quadratic susceptibility.

\subsection{Numerical results}
We use the same grids, frequencies, sources, noise level and relative cutoff
($0.03$) as in Section~\ref{sec:numerics}, now adding noise to both the
fundamental and the second harmonic data. The nonlinear data are computed on
the fine grid by iterating on the quadratic terms, with the linear
susceptibility treated implicitly. The linear susceptibility is
\[
\eta_1(x,y)=A\exp\left(-\frac{(x-0.62)^2+(y-0.62)^2}{2(0.08)^2}\right),
\]
and $\eta_2$ is $0.1$ times the Gaussian \eqref{eq:gaussian-target}. Since the
fundamental equation can be written as
$(\Delta+k^2(1+4\pi\eta_1))u_1=-8\pi k^2\eta_2u_2\overline{u_1}$, the relative
permittivity is $1+4\pi\eta_1$, so $A=0.02$ and $A=0.06$ correspond to a peak
relative permittivity of $1.25$ and $1.75$. The truncated SVD
keeps $113$ singular values for $\eta_1$ and $277$ for $\eta_2$. Fewer are
kept for $\eta_1$, because the fundamental data are at the lower frequency
$k$, and the projection error is correspondingly larger for $\eta_1$ than
for $\eta_2$.

\begin{table}[ht]
\centering
\begin{tabular}{lrrrrrr}
\toprule
 & IBS 1 & IBS 2 & IBS 3 & IBS 4 & IBS 5 & IBS 7\\
\midrule
\multicolumn{7}{l}{$A=0.02$, $\gamma=1$}\\
\quad$\eta_1$, rel.\ to projection & $6.2\%$ & $0.9\%$ & $1.3\%$ & $1.3\%$ & $1.4\%$ & $1.4\%$\\
\quad$\eta_2$, rel.\ to projection & $9.6\%$ & $2.3\%$ & $1.4\%$ & $1.4\%$ & $1.4\%$ & $1.4\%$\\
\multicolumn{7}{l}{$A=0.06$, $\gamma=1$}\\
\quad$\eta_1$, rel.\ to projection & $18.9\%$ & $6.3\%$ & $1.8\%$ & $2.0\%$ & $3.5\%$ & $3.7\%$\\
\quad$\eta_2$, rel.\ to projection & $29.0\%$ & $11.4\%$ & $5.1\%$ & $4.7\%$ & $3.9\%$ & $3.5\%$\\
\multicolumn{7}{l}{$A=0.06$, $\gamma=0.5$}\\
\quad$\eta_1$, rel.\ to projection & $18.8\%$ & $6.3\%$ & $1.8\%$ & $2.0\%$ & $3.5\%$ & $3.7\%$\\
\quad$\eta_2$, rel.\ to projection & $29.2\%$ & $11.5\%$ & $5.0\%$ & $4.7\%$ & $3.9\%$ & $3.5\%$\\
\bottomrule
\end{tabular}
\caption{Full SHG model: relative $L^2$ error of the partial sums of the IBS,
measured relative to the projections of $\eta_1$ and $\eta_2$ onto the
retained singular vectors.}
\label{tab:full}
\end{table}

\begin{figure}[ht]
\centering
\includegraphics[width=\textwidth]{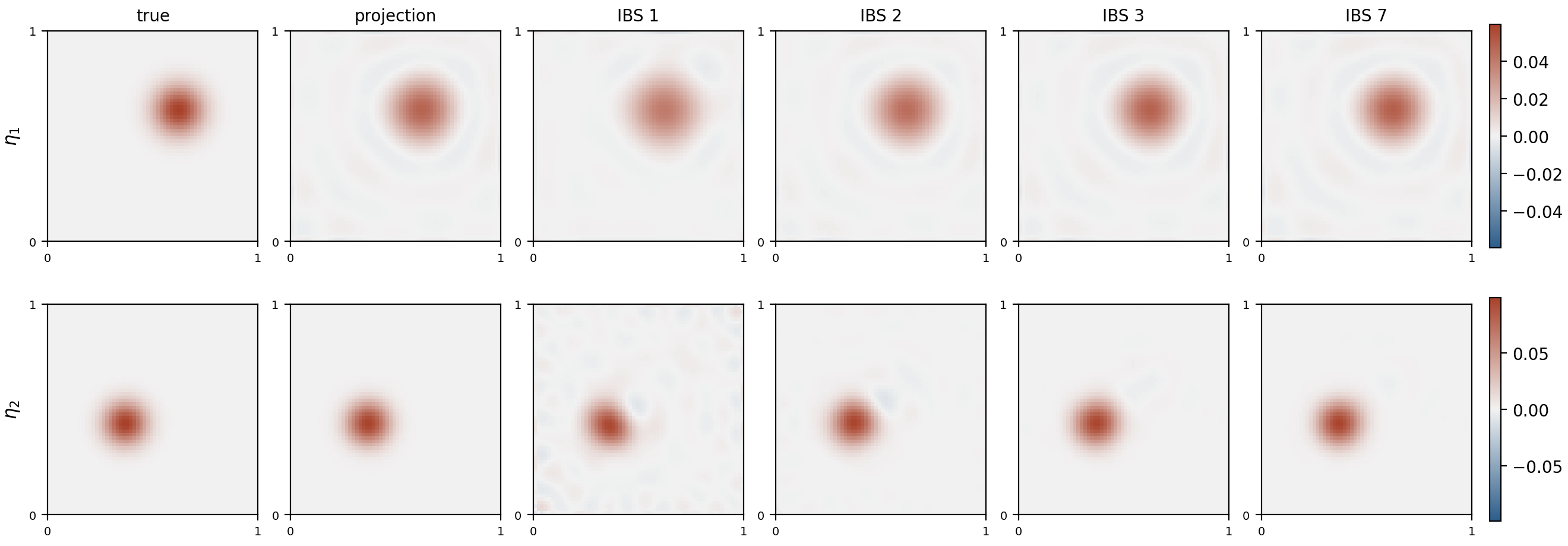}
\caption{Full SHG model with $A=0.06$ (a peak relative permittivity of $1.75$) and
$\gamma=1$: reconstruction of the
linear susceptibility $\eta_1$ (top) and the quadratic susceptibility $\eta_2$
(bottom). The data are computed on a $128\times128$ grid and the inversion
uses a $64\times64$ grid.}
\label{fig:full}
\end{figure}

The results are shown in Table~\ref{tab:full} and Figure~\ref{fig:full}. With
$A=0.02$ and $\gamma=1$, the nonlinear part is about one eighth of the data
at both frequencies. In the first term, the reconstruction of $\eta_2$
contains artifacts caused by the coupling between the two susceptibilities,
since the second harmonic data depend on $\eta_1$ at second order. The higher
order terms remove these errors, and the series converges to within the
level of the model error and noise of the projections of both
susceptibilities. Thus the two susceptibilities can be recovered
simultaneously, and the IBS corrects for their coupling.

With a larger linear susceptibility, $A=0.06$ and $\gamma=1$, the nonlinear
part is about a third of the data, and the higher order terms are essential
(see Figure~\ref{fig:full}). The first term underestimates $\eta_1$ and
smears it, and its reconstruction of $\eta_2$ is distorted, with an artifact
on the side facing $\eta_1$. The partial sum of order three reduces the
errors relative to the projections several times over, and the artifact is
removed. After that the error for $\eta_2$ continues to decrease slowly,
while the error for $\eta_1$ increases slowly, so the series is near the
edge of convergence.

Finally, we repeat this experiment with the weaker source $\gamma=0.5$. In
the reduced model this would halve the effective contrast. Here the results
are essentially unchanged (Table~\ref{tab:full}). This is what we expected in
Step~3: the nonlinear effects in this experiment are dominated by the linear
susceptibility and its coupling with $\eta_2$, and these do not become weaker
when the source is made weaker. Convergence for the full model therefore
requires a bound on the linear susceptibility, while the quadratic
susceptibility can still be controlled by the source strength.

A theory for the full model, with bounds on the real-multilinear forward
operators and convergence conditions that separate the linear and quadratic
contrasts, is left for future work.

\section*{Code availability}
The Python code and the scripts that reproduce the numerical results will be
made publicly available at
\begin{center}
\url{https://github.com/sharimoskow/shg-ibs}.
\end{center}
The script \texttt{run\_fine\_2d.py} produces the results in
Section~\ref{sec:numerics}, the script \texttt{run\_full\_shg\_2d.py} those in
Section~\ref{sec:full}, and \texttt{make\_paper\_figures.py} the figures.

\section*{Acknowledgments}
S. Moskow was partially supported by NSF grant DMS-2308200.
J. C. Schotland was partially supported by the Simons Foundation grant SFI-MPS-TSRG-00023182-04.

Generative AI systems (OpenAI ChatGPT/Codex and Anthropic Claude) were used
as assistants in the preparation of this manuscript: for drafting and editing
the text, for writing and debugging the numerical code, and for checking and
refining notation, calculations and proofs. The authors verified all of the
mathematical arguments, numerical results and text, and take full
responsibility for the contents.

\end{document}